\documentclass{amsproc}

\usepackage{amsmath,amssymb,amscd,color}
\usepackage{graphicx,epstopdf}
\usepackage{float}
\graphicspath{ {./} }

\def\e{\epsilon}

\def\bp{\begin{proposition}}
\def\ep{\end{proposition}}
\def\bt{\begin{theo}}
\def\et{\end{theo}}
\def\be{\begin{equation}}
\def\ee{\end{equation}}
\def\bl{\begin{lemma}}
\def\el{\end{lemma}}
\def\bc{\begin{corollary}}
\def\ec{\end{corollary}}
\def\pr{\noindent{\bf Proof: }}

\def\bd{\begin{definition}}
\def\ed{\end{definition}}

\newcommand{\x}[1]{{}$\kern-2\mathsurround${}\binoppenalty10000 \relpenalty10000 #1{}$\kern-2\mathsurround${}}

\makeatletter
\DeclareRobustCommand*\cal{\@fontswitch\relax\mathcal}
\makeatother

\newtheorem{theo}{Theorem}[section]
\newtheorem{lemma}{Lemma}[section]
\newtheorem{definition}{Definition}[section]
\newtheorem{corollary}{Corollary}[section]
\newtheorem{proposition}{Proposition}[section]

\numberwithin{equation}{section}

\begin{document}

\title[Accuracy of reconstruction of spike-trains]{Accuracy of reconstruction of spike-trains with two near-colliding nodes}

\author[A. Akinshin]{Andrey Akinshin}
\address{Department of Mathematics,	The Weizmann Institute of Science, Rehovot 76100, Israel}
\address{Laboratory of Inverse Problems of Mathematical Physics, Sobolev Institute of Mathematics SB RAS, Novosibirsk 630090, Russia}
\email{andrey.akinshin@weizmann.ac.il}
\thanks{}

\author[G. Goldman]{Gil Goldman}
\address{Department of Mathematics,	The Weizmann Institute of Science, Rehovot 76100, Israel}
\email{gilgoldm@gmail.com}
\thanks{}

\author[V. Golubyatnikov]{Vladimir Golubyatnikov}
\address{Laboratory of Inverse Problems of Mathematical Physics, Sobolev Institute of Mathematics SB RAS, Novosibirsk 630090, Russia}
\email{golubyatn@yandex.ru}
\thanks{}

\author[Y. Yomdin]{Yosef Yomdin}
\address{Department of Mathematics,	The Weizmann Institute of Science, Rehovot 76100, Israel}
\email{yosef.yomdin@weizmann.ac.il}

\thanks{The work was supported by the RFBR grant 15-01-00745 A; ISF, Grant No. 779/13.}

\date{}

\begin{abstract}
We consider a signal reconstruction problem for signals $F$ of the form $ F(x)=\sum_{j=1}^{d}a_{j}\delta\left(x-x_{j}\right),$
from their moments $m_k(F)=\int x^kF(x)dx.$ We assume $m_k(F)$ to be known for $k=0,1,\ldots,N,$ with an absolute error not exceeding $\e > 0$.

We study the ``geometry of error amplification'' in reconstruction of $F$ from $m_k(F),$
in situations where two neighboring nodes $x_i$ and $x_{i+1}$ near-collide, i.e $x_{i+1}-x_i=h \ll 1$.
We show that the error amplification is governed by certain algebraic curves $S_{F,i},$ in the parameter space of signals $F$,
along which the first three moments $m_0,m_1,m_2$ remain constant.
\end{abstract}

\maketitle

\section{Introduction}\label{Sec:Intro}


The problem of reconstruction of spike-trains, and of similar signals, from noisy moment measurements,
and a closely related problem of robust solving the classical Prony system,
is a well-known problem in Mathematics and Engineering.
It is of major practical importance, and, in case when the nodes nearly collide,
it presents major mathematical difficulties.
It is closely related to a spike-train ``super-resolution problem'', (see \cite{Aza.Cas.Gam,Bat1,Can.Fer1,Can.Fer2,Con.Hir,Dem.Ngu,Dem.Nee.Ngu,Don,Duv.Pey,Fer,Hec.Mor.Sol,Lia.Fan,Mor.Can,Moi,Pet.Plo} as a small sample).

The aim of the present paper is to investigate the possible amplification of the measurements error $\e$ in the
reconstruction process, caused by the fact that some of the nodes of $F$ near-collide. Recently this problem attracted
attention of many researchers. In particular, in \cite{Aki.Bat.Yom,Bat1,Can.Fer2,Dem.Nee.Ngu,Don} it was shown (in
different settings of the problem) that if $s$ spikes of $F$ are near-colliding in an interval of size $h\ll 1$, then a strong ``noise
amplification'' occurs: up to a factor of $(\frac{1}{h})^{2s-1}$.
Specifically, in \cite{Aki.Bat.Yom} a parametric setting (the same as in the present paper) was considered (see, as a
small sample, \cite{Bat.Yom3,Bey.Mon2,Blu.Dra.Vet.Mar.Cou,Pet.Plo,Pet.Pot.Tas,Plo.Wis,Yom} and references therein). 
In this setting, signals are assumed to be members of a parametric family with a finite number of parameters.
The parameters of the signal are then considered as unknowns, while the measurements provide a system of algebraic equations in these unknowns.

It was announced in \cite{Aki.Bat.Yom} that the strongest ``noise amplification'' 
occurs along the algebraic curves $S$ (``Prony curves''), defined in the signal parameter space by the $2s-1$ initial equations of the classical ``Prony system'' (system \ref{eq:Prony.system1} below).
However, \cite{Aki.Bat.Yom} provided neither detailed proofs, nor explicit constants in the error bound, nor the explicit description of the curves $S$.

In the present paper we consider reconstruction of spike-train signals of an a priori known form
$F(x)=\sum_{j=1}^{d}a_{j}\delta\left(x-x_{j}\right)$, from their moments $m_0(F),\ldots,\allowbreak m_{N-1}(F)$, $ N\ge 2d$, in the
case where two nodes $x_i,x_{i+1}$ near-collide. That is, $x_{i+1}-x_{i+1}=h\ll 1$.

\smallskip

In Section \ref{Sec:setting} we introduce the $\e$-error set $E_\e(F)$, consisting of all signals $F'$, for which the moments of $F'$ differ from the moments of $F$ by at most $\e$. 
The set $E_\e(F)$ presents the distribution of all the possible reconstructed signals $F'$, 
caused by independent errors, not exceeding $\e$, in each of the moment measurements $m_k(F)$. 
Thus the geometry of $E_\e(F)$ reflects the patterns of the possible error amplification in the reconstruction process.

In this paper we are mostly interested in the lower bounds for the error in nodes reconstruction. 
We thus consider the projection $E^x_\e(F)$ of the error set $E_\e(F)$ to the nodes space. 
This set represents the error amplification in nodes reconstruction. 
In particular, the ``radius'' $\rho^x_\e(F)$ of $E^x_\e(F)$ provides a lower bound on the nodes reconstruction accuracy of any reconstruction algorithm 
(see a more detailed description of this fact in Section \ref{Sec:setting}).

\smallskip

One of our two main results - Theorem \ref{thm:main2} in Section \ref{Sec:Error.amplification.main} - is that {\it for $F$ with two nodes in a distance $h$, and for any $\e$ of order $h^3$ (or larger) we have $\rho^x_\e(F)\ge Ch.$ Consequently, the presence of near-colliding nodes implies a massive amplification of the measurements error in the process of nodes reconstruction - up to $h^{-2}$ times.}

\smallskip

In order to prove Theorem \ref{thm:main2} we start in this paper the investigation of the geometry of the error sets $E_\e(F)$ and $E^x_\e(F)$ (which, as we believe, is important by itself). 
First we provide in Section \ref{sec:examples} numerical simulations and visualizations, 
which suggest that for $\e \sim h^3$ the $\e$-error set $E_\e(F)$ is an ``elongated curvilinear parallelepiped'' of the width $\sim h$, 
stretched up to the size $\sim 1$ along a certain curve $S_F$ (while its projection onto the nodes space, $E^x_\e(F)$, is stretched along the projection $S^x_F$ of $S_F$). 
These experiments suggest also that as $h\to 0$, the sets $E_\e(F)$ and $E^x_\e(F)$ concentrate closer and closer around the curves $S_F$ (respectively, $S^x_F$).

\smallskip

Next, we give in Section \ref{Sec:Error.amplification.main} an independent definition of the ``Prony curves'' $S_F$, ``discovered'' in Section \ref{sec:examples}: for each $F$ the Prony curve $S_F,$ passing through $F$ in the signal parameter space, is defined by the requirement that along it the first three moments $m_0,m_1,m_2$ do not change. An explicit parametric description of the curves $S_F$ is given in Section \ref{sec:parm.Prony.curves}.

\smallskip

Our second main result - Theorem \ref{thm:main1} in Section \ref{Sec:Error.amplification.main} - is that {\it indeed, as suggested by visualizations in Section \ref{sec:examples}, the set $E_\e(F)$ contains a ``sufficiently long'' part of the Prony curve $S_F$ around $F$}. As a consequence, we obtain Theorem \ref{thm:main2}. Let us stress that all the constants in Theorems \ref{thm:main1} and \ref{thm:main2} are explicit (and reasonably realistic).

\smallskip

The proofs are given in Section \ref{Sec:proofs}.

\smallskip

Finally, in Section \ref{sec:compl.nodes}, we compare two approaches to the reconstruction problem for real spike-train signals:
from their moments, and from their Fourier samples (which can be interpreted as the moments of an
appropriate signal $\tilde F$ with complex nodes).
Recently in \cite{Aki.Gol.Yom} a trigonometric reconstruction method for $\tilde F$ was suggested,
which uses, as an input, {\it only three complex moments
$m_0(\tilde F),m_1(\tilde F),m_2(\tilde F)$}.
According to the approach of the present paper, we would expect for the trigonometric method
(for $\tilde F$ with two nodes in a distance $h\ll 1$, and for $\e\sim h^3$,) the worst
case reconstruction error of order $\sqrt \e$, while for the Prony inversion we show it to be of order $\e^{\frac{1}{3}}$.
We pose some open questions related to this apparent contradiction.

\medskip

The authors would like to thank the referee for a constructive criticism,
as well as for remarks and suggestions, which allowed us to significantly improve the presentation.


\section{Setting of the problem}\label{Sec:setting}

Assume that our signal $F(x)$ is a spike-train, that is, a linear combination of $d$ shifted $\delta$-functions:
\be \label{eq:equation.model.delta}
F(x)=\sum_{i=1}^{d}a_{i}\delta\left(x-x_{i}\right),
\ee
where $a=(a_1,\ldots,a_d) \in {\mathbb R}^d, \ x=(x_1,\ldots,x_d) \in {\mathbb R}^d.$
We assume that the form (\ref{eq:equation.model.delta}) is a priori known, but the specific parameters $(a,x)$ are unknown.
Our goal is to reconstruct $(a,x)$ from $N\ge 2d$ moments $m_k(F)=\int_{-\infty}^\infty x^k F(x)dx, \ k=0,\ldots,N-1$, which are known with a possible absolute error 
of no more than $\e>0$.

\smallskip

The moments $m_k(F)$ are expressed through the unknown parameters $(a,x)$ as $m_k(F)=\sum_{i=1}^d a_i x_i^k$.
Hence our reconstruction problem is equivalent to solving the (possibly over-determined) {\it Prony system} of algebraic equations, with the unknowns $a_i,x_i$:

\be\label{eq:Prony.system1}
\sum_{i=1}^d a_i x_i^k = m_k(F), \ k= 0,1,\ldots,N-1.
\ee
This system appears in many theoretical and applied problems.
There exists a vast literature on Prony and similar systems,
in particular, on their robust solution in the presence of noise - see, as a small sample,
\cite{Aza.Cas.Gam,Bey.Mon2,Lia.Fan}, \cite{Ode.Bar.Pis}-\cite{Pot.Tas}, and references therein.

We present the spike train reconstruction problem in a geometric language of spaces and mappings.
Let us denote by ${\cal P}={\cal P}_d$ the parameter space of signals $F$,
$$
{\cal P}_d=\{(a,x)=(a_1,\ldots,a_d,x_1,\ldots,x_d)\in {\mathbb R}^{2d}, \ x_1<x_2<\ldots<x_d \},
$$
and by ${\cal M}={\cal M}_{N} \cong {\mathbb R}^{N}$ the moment space, consisting of the $N$-tuples of moments $(m_0,m_1,\ldots,m_{N-1})$.
We will identify signals $F$ with their parameters $(a,x)\in {\cal P}.$

\smallskip

The Prony mapping $PM=PM_{d,N}:{\cal P}_d\to {\cal M}_{N}$ is given by
$$
PM(F)= \mu = (\mu_0,\ldots,\mu_{N-1}) \in {\cal M}, \ \mu_k=m_k(F), \ k=0,\ldots,N-1.
$$
Inversion of the Prony mapping is equivalent to reconstruction of a spike-train signal $F$ from its moments (or to solving Prony system (\ref{eq:Prony.system1})).

\smallskip

The aim of this paper is to investigate the amplification of the measurements error $\e$ in the reconstruction process, in case of two
near-colliding nodes.
We are interested in effects, caused by the geometric nature of system (\ref{eq:Prony.system1}), independently of the specific method of its inversion.

\smallskip

The error amplification is reflected by the geometry of the $\e$-error set $E_\e(F)$, which is defined as follows:

\bd\label{def:error.set}
The $\e$-error set $E_\e(F)$ consists of all signals $F'\in {\cal P}_d$, for which the moments of $F'$ differ from the moments of $F$ by at most $\e$:
$$
E_\e(F) = \{F'\in {\cal P}_d, \ |m_k(F')-m_k(F)|\le \e, \ k=0,\ldots, N-1\}.
$$
Equivalently, $E_\e(F)=PM_{d}^{-1}(Q^{N}_\e(F))$, where $Q^{N}_\e(F)\subset {\cal M}_{N}$ is the $N$-dimensional $\e$-cube centered at $PM(F)\in {\cal M}_{N}$.
\ed
The $\e$-error set $E_\e(F)$ presents the distribution of possible reconstructed signals $F'$,
caused by the independent errors, not exceeding $\e$, in each of the moment measurements $m_k(F)$.
Its yet another convenient description is as the set of solutions of the Prony system
\be\label{eq:Prony.system2}
\sum_{i=1}^d a_i x_i^k = m_k(F)+\epsilon_k, \ k= 0,1,\ldots,N-1,
\ee
with all the possible errors $\epsilon_k$ satisfying $|\epsilon_k|\leq \e, \ k= 0,1,\ldots,N-1.$

Notice that the $\e$-error set $E_\e(F)$ depends on $N$, the number of the moments which we use in reconstruction.
Since $N$ is assumed to be fixed, we do not indicate it in the notations.
\smallskip

In this paper we are mainly interested in the accuracy of the nodes reconstruction which is determined by the geometry of the projection $E^x_{\epsilon}(F)$ of the set $E_{\epsilon}(F)$ onto the nodes space. Accordingly, we define the worst case error $\rho^x_\e(F)$ in the Prony reconstruction of the nodes of $F$ as follows:

\bd\label{def:wrst.case.error}
For $F=(a,x)\in {\cal P}_d$, the worst case error $\rho^x_\e(F)$ in the reconstruction of the nodes of $F$ is defined by
\be\label{eq:worst.case.error}
\rho^x_\e(F)=\sup_{F'=(a',x')\in E_\e(F)} ||x'-x||,
\ee
where $||\cdot||$ denotes the Euclidean norm in the space of the nodes.
\ed
In fact, $\frac{1}{2}\rho^x_\e(F)$ bounds from below the worst case error in nodes reconstruction with {\it any} reconstruction algorithm $A$. Indeed, we can informally argue as follows: let $F''=(a'',x'')\in E_\e(F)$ be a signal for which the supremum in (\ref{eq:worst.case.error}) is nearly achieved. Assume that we apply $A$ to both signals $F$ and $F''$, and the (adversary) noise is zero for $F$ and is equal to the difference of the moments of $F$ and $F''$ in the second case. Thus $A$ obtains as the input in both cases the moments of $F$. Whatever result $\tilde F=(\tilde a, \tilde x)$ the algorithm $A$ produces as an output, at least one of the distances $||x-\tilde x||$ or $||x''-\tilde x||$ will be not smaller than $\frac{1}{2}\rho^x_\e(F)$.

\section{Visualization of the error sets  $E_\epsilon(F)$}\label{sec:examples}

For a given $h, \ 0<h<1,$ consider a signal $F(x)=\frac{1}{2}\delta(x+h)+\frac{1}{2}\delta(x-h) \in {\cal P}_2.$ We put $N=2d=4$. The moments of $F$ are
$$
m_0(F)=1, \ m_1(F)=0, \ m_2(F)=h^2, \ m_3(F)=0.
$$
For a given $\e>0$ we consider the $\e$-cube $Q^4_\e(F)\subset {\cal M}_{4}$ centered at $(1,0,h^2,0)\in {\cal M}_{4}$, 
and the $\e$-error set $E_\e(F)=PM_{2,4}^{-1}(Q^4_\e(F)).$ Equivalently, $E_\e(F)$ is defined in ${\cal P}_2$ by the inequalities
$$
|m_0(F')-1| \le \e, \ |m_1(F')| \le \e, \ |m_2(F')-h^2| \le \e, \ |m_3(F')| \le \e.
$$

\begin{figure}
		\centering
		\includegraphics[scale=0.80]{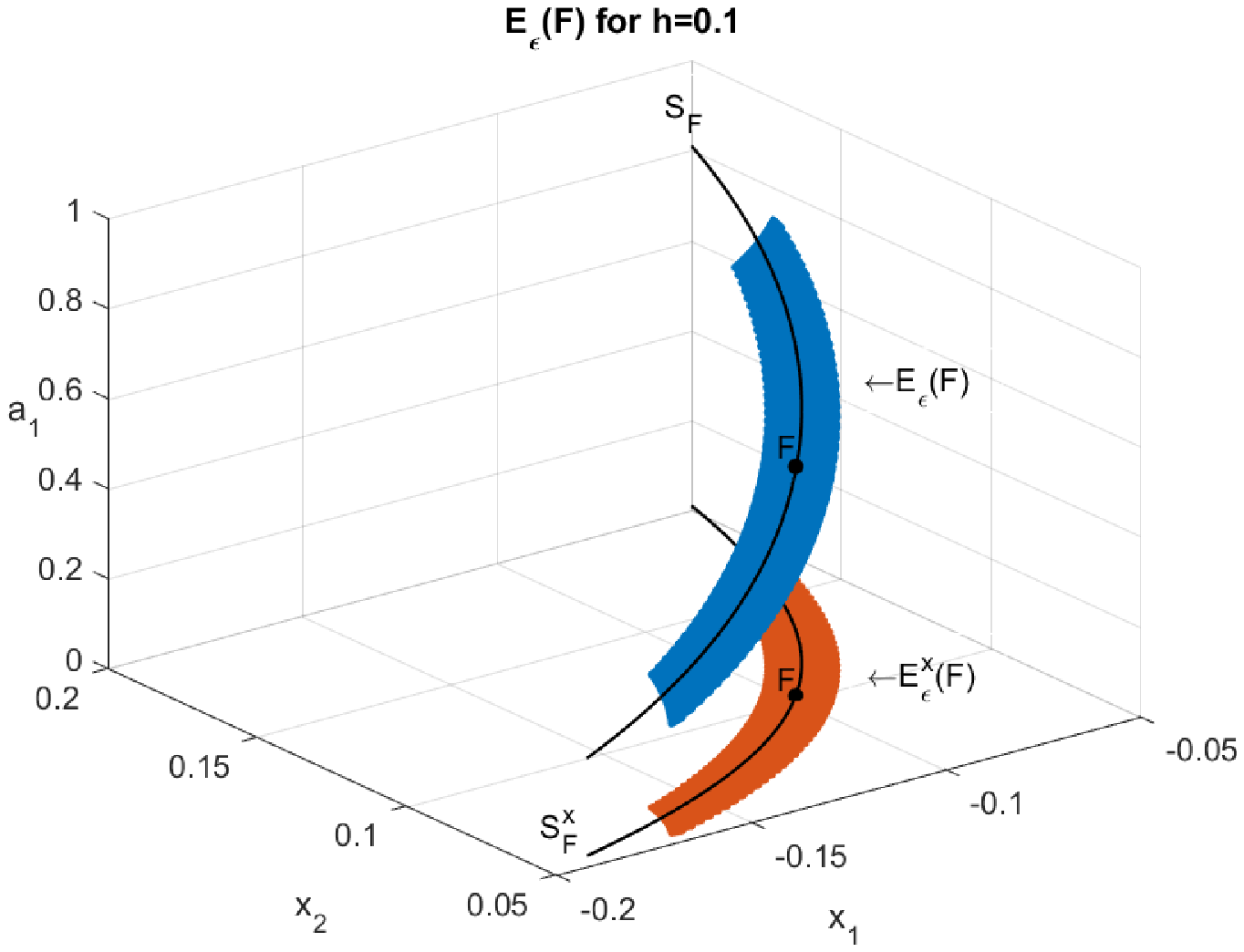}
		\caption{The error set $E_{\e}(F)$ and its projection $E_{\e}^x(F)$ for
		\;\;\;\;\;$h=0.1$, $\e=2h^3=0.002$ and
		$F(x)=\frac{1}{2}\delta(x-0.1)+\frac{1}{2}\delta(x+0.1)$.}
		\label{fig:figure1}
		\hspace{0.5cm}
		\centering
		\includegraphics[scale=0.80]{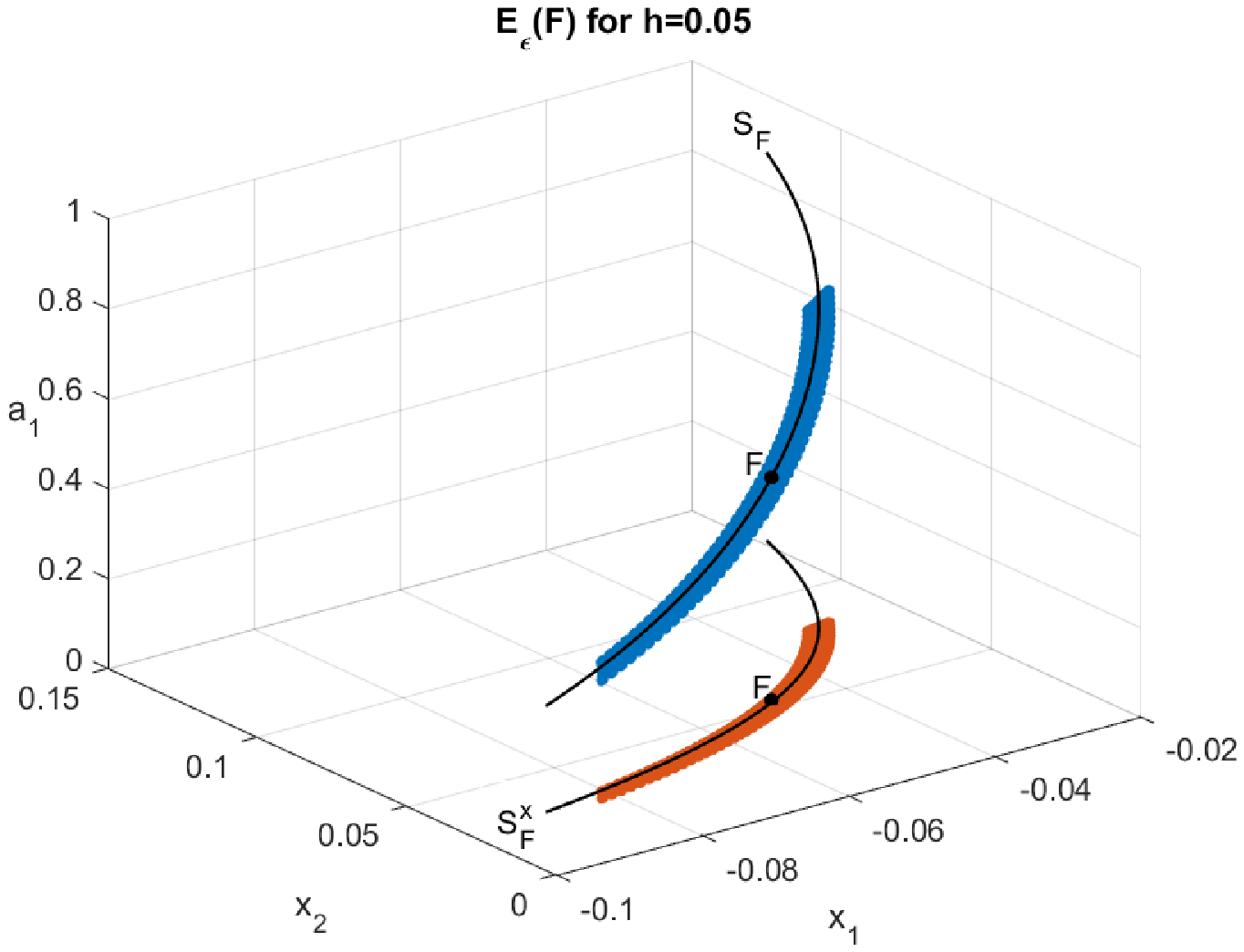}
		\caption{The error set $E_{\e}(F)$ and its projection $E_{\e}^x(F)$ for
		\;\;\;\;\;$h=0.05$, $\e=2h^3=0.00025$ and
		$F(x)=\frac{1}{2}\delta(x-0.05)+\frac{1}{2}\delta(x+0.05)$.}
		\label{fig:figure2}
\end{figure}

The $\e$-error set $E_\e(F)$ is a four-dimensional subset of ${\cal P}_2\cong
{\mathbb R}^4,$ and its direct visualization is problematic. Instead the
following Figures \ref{fig:figure1} and \ref{fig:figure2} show the projection of
$E_\e(F)$ onto the three-dimensional coordinate subspace of ${\cal P}_2$, spanned by the two nodes coordinates $x_1,x_2$ and the first amplitude $a_1$, as well as its further projection $E^x_\e(F)$ onto the nodes plane $x_1,x_2$.

Notice, that by the first of the Prony equations $a_1+a_2=m_0(F)+\e_0$ in (\ref{eq:Prony.system2}) we have $a_2=m_0(F)-a_1+\e_0,$ with $|\e_0|\le \e.$ Thus the projections of $E_\e(F)$ shown in Figures \ref{fig:figure1} and \ref{fig:figure2}, give a rather accurate (up to $\e$) representation of the true error set.

\medskip


\medskip



Let us stress a natural scaling in our problem, reflected in Figures \ref{fig:figure1} and \ref{fig:figure2}: {\it the scale in nodes is of order $h$, while the scale in the amplitudes is of order $1$}.

\medskip

Figures \ref{fig:figure1} and \ref{fig:figure2} suggest that $E_\e(F)$ is an
``elongated curvilinear parallelepiped'', with the sizes of its two largest edges of orders $1$ and $h,$ respectively. (The third and the fourth edges, of orders $h^2$ and $h^3,$ respectively, are not visible). $E_\e(F)$ is stretched up to the size $\sim 1$ along a certain curve $S$, depicted in the pictures. Respectively, $E^x_\e(F)$ is stretched up to the size $\sim h$ along the projection curve $S^x$. A comparison between Figures \ref{fig:figure1} and \ref{fig:figure2} also suggests that as $h$ (and $\e \sim h^3)$ decrease, the error set concentrates closer and closer along the curve $S$. (Compare a conjectured general description of $E_\e(F)$ at the end of Section \ref{Sec:Error.amplification.main}).

\smallskip

Below we analyse the structure of the error set $E_\e(F)$ in some detail, and show that {\it $S$ is an algebraic curve, which we call the ``Prony curve''}. We show that for $F$ as above, the projection $S^x$ of $S$ onto the node subspace is the hyperbola $x_1x_2=-h^2$, while $a_1$ is expressed on this curve through $x_1,x_2$ as $a_1=\frac{x_2}{x_2-x_1}.$ In Section \ref{sec:parm.Prony.curves} we study such ``Prony curves'' in detail.

\smallskip

Numerically, the figures above were constructed via the following procedure:
we construct a four-dimensional regular net $Z\subset Q^4_\e(F)\subset {\cal M}_{2},$ with a sufficiently small step.
For each point $z\in Z,$ its Prony preimage $w=PM^{-1}(z)\in {\cal P}_2$ is calculated, and the projection of $w$ onto the space $(a_1,x_1,x_2)$ is plotted.

Some other visualisation results can be found in \cite{Aki}.

\smallskip

In what follows we assume that {\it inversion of the Prony map, or solving of (\ref{eq:Prony.system2})
(when possible) is accurate, and the reconstruction error is caused only by the measurements errors $\epsilon_k$.}

\section{Prony curves and error amplification: main results}\label{Sec:Error.amplification.main}

We will consider signals $F(x)$ of the form (\ref{eq:equation.model.delta}),
with two near colliding nodes $x_i$ and $x_{i+1}, \ 1 \leq i\leq d-1.$
In the present paper we study the geometry of the reconstruction error,
allowing perturbations only of the cluster nodes $x_i,x_{i+1},$ and of their amplitudes $a_i,a_{i+1}$.
Therefore the positions and the amplitudes of the other nodes are not relevant for our results.
However, in order to avoid possible collisions of the cluster nodes with their neighbors in the process of deformation,
we will always assume that for $x_{i+1}-x_i=h>0$, the distances to the neighboring nodes from the left and from the right satisfy
$x_i-x_{i-1}\geq 3h, \ x_{i+1}-x_{i+2}\geq 3h.$
We do not assume formally that $h \ll 1,$ but this is the case where the geometric patterns we describe become apparent.

\smallskip

For each signal $F$ and index $i$ the Prony curve $S=S_{F,i}$ passing through $F$ is obtained by varying only the nodes and amplitudes $(a_i,x_i),(a_{i+1},x_{i+1})$,
while preserving the first three moments. More accurately, we have the following definition:

\bd\label{def:Prony.curves}
Let $F(x)=\sum_{j=1}^{d}a_{j}\delta\left(x-x_{j}\right)\in {\cal P}_d$ and let $i, \ 1 \leq i\leq d-1,$ be fixed.
The Prony curve $S=S_{F,i}\subset {\cal P}_d$ consists of all the signals
$$
F'(x)= \sum_{j=1}^{d}a'_{j} \delta(x-x'_{j})\in {\cal P}_d
$$
for which $a'_j=a_j, x'_j=x_j, \ j\ne i,i+1,$ and $m_k(F')=m_k(F)$ for $k=0,1,2.$
\ed
By definition, we always have $F\in S_{F,i}$.
In this paper we concentrate on an ``$h$-local'' part $S_{F,i}(h)$ around $F$ of the Prony curve $S_{F,i}$, consisting of all
$$
F'(x)= \sum_{j=1}^{d}a'_{j} \delta(x-x'_{j})\in S_{F,i},
$$
for which the nodes $(x'_i,x'_{i+1})$ belong to the disk $D \subset \mathbb{R}^2$ of radius $\frac{1}{2}h$ centered at $(x_i,x_{i+1}).$
In particular, the node collision cannot happen on $S_{F,i}(h)$ (see Lemma \ref{lem:bdd.dist11} below). Compare also to Figure \ref{fig:figure3}.

\smallskip

Notice, however, that the Prony curves $S_{F,i}$ are {\it global algebraic curves (possibly singular)}.
Their explicit global parametrization is described in Section \ref{sec:parm.Prony.curves} below, and one can show that in some cases they can pass through the node collision points (with amplitudes tending to infinity).
We believe that the Prony curves and their multi-dimensional generalizations play an important role in understanding of multi-nodes collision singularities.


\medskip

The following definition specifies the type of signals we will work with:

\bd\label{def:i.h.M.clusters}
A signal $F=\sum_{j=1}^{d}a_{j}\delta\left(x-x_{j}\right)\in {\cal P}_d$ is said to form an $(i,h,M)$-cluster,
with given $h, \ 0< h <1, \ M>0,$ and $i, \ 1 \leq i\leq d-1,$ if $x_{i+1}-x_i=h,$ and $|a_i|,|a_{i+1}|\leq M$.
\ed
For $F$ forming an $(i,h,M)$-cluster, let $\kappa=\frac{x_{i+1}+x_i}{2}$ be the center of the interval $[x_i,x_{i+1}]$.
We put $C(F)=18M(1+|\kappa|)^N$ (Notice that $C(F)$ depends on $\kappa,M$, but not on $h$). One of our main results is the following:

\bt\label{thm:main1}
Let $F\in {\cal P}_d$ form an $(i,h,M)$-cluster. Then for each $\e\ge C(F)h^3$ the $\e$-error set $E_\e(F)$ contains the local Prony curve $S_{F,i}(h).$
\et
\begin{figure}[H]
			\centering
			\includegraphics[width=\textwidth]{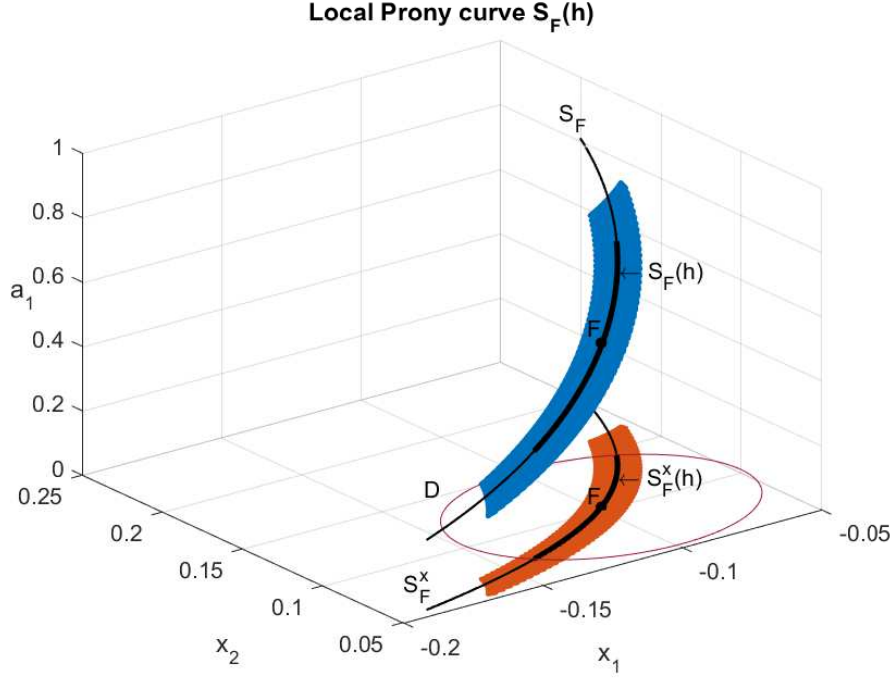}
			\caption{The error set $E_{\e}(F)$ and its projection $E_{\e}^x(F)$ for
		\;\;\;\;\;$h=0.1$, $\e=2h^3=0.002$ and
		$F(x)=\frac{1}{2}\delta(x-0.1)+\frac{1}{2}\delta(x+0.1)$. The circle depicts
		the boundary of the disk $D$. In bold is the local prony curve $S_{F,i}(h)$ and
		its projection.}
			\label{fig:figure3}
\end{figure}
Figure \ref{fig:figure3} suggests that in fact the $\e$-error set $E_\e(F)$ ``concentrates'' around the local Prony curve $S_{F,i}(h).$
Already the fact that this curve is inside $E_\e(F)$ (provided by Theorem \ref{thm:main1}) implies important conclusion on the worst case reconstruction error.
Indeed, as we show below, the projection $S^x_{F,i}(h)$ of the local Prony curve $S_{F,i}(h)$ onto the nodes space has a length of order $h$.
Consequently, in the presence of an $h$-cluster in $F$, and for each $\e\ge C(F)h^3$,
there are signals $F'\in E_\e(F),$ with the nodes $x'_i,x'_{i+1}$ at a distance $\sim h$ from $x_i,x_{i+1}$.
That is, a massive error amplification from $h^3$ to h occurs in this case.

\smallskip




Our second main result presents this fact formally, in terms of the worst case error:

\bt\label{thm:main2}
Let $F\in {\cal P}_d$ form an $(i,h,M)$-cluster.
Then for each $\e\ge C(F)h^3$ the worst case reconstruction error $\rho^x_\epsilon(F)$ in nodes of $F$ is at least $\frac{1}{2}h$.
\et
We see that a measurements error $\e\sim h^3$ can be amplified up to the factor $\sim h^{-2}$ in reconstruction of the nodes of $F$.
In particular, for $d=2,$ i.e, in the case of exactly two nodes in $F$,
and for $N=4$, we get (assuming that $M=1$ and $x_1,x_2\in [-1,1])$ that $C(F)\le288.$ So the minimal accuracy $\e$
required to keep the error in the nodes reconstruction less than $\frac{1}{2}h$ is $\e \le 288 h^3$.
For $h=0.01$ we get $\e\leq 0.0003$.

\smallskip

In \cite{Aki.Bat.Yom}, in order to show that the length of the curve $S^x_{F,i}(h)$ is of order $h$,
we use the inverse function theorem, combined with estimations from \cite{Bat.Yom2} of the Jacobian of the Prony mapping.
As a result, the constants become much less explicit.

\smallskip

We expect that the results of Theorems \ref{thm:main1} and \ref{thm:main2} can be extended
to an accurate description of the $\e$-error set in the case of clusters with more than two nodes,
using an appropriate version of the ``quantitative inverse function theorem''. Informally, we expect the following general result to be true:

\smallskip

\noindent {\it Let the nodes $x_1,\ldots,x_d$ of $F$ form a cluster of size $h\ll 1.$ Then for $\e \le O(h^{2d-1})$ the $\e$-error set $E_\e(F)$ is a ``non-linear coordinate parallelepiped'' $\Pi_{h,\e}(F)$ with respect to the moment coordinates $m_k(F')$, centered at $F$. Its width in the direction of the moment coordinate $m_k, \ k=0,\ldots,2d-1,$ is of order $\e h^{-k}.$ In particular, the maximal stretching of $\Pi_{h,\e}(F)$, of order $\e h^{-(2d-1)},$ occurs along the Prony curve $S_{2d-2}(F)$.}

\smallskip

However, an application of the approach based on the inverse function theorem will significantly reduce the domain of applicability of the results, and will make the constants less explicit.

On the other hand, we believe that the explicit parametric description of the Prony curves, given in Section \ref{sec:parm.Prony.curves} below, can be extended to clusters with more than two nodes. It becomes significantly more complicated, but promises potentially better understanding of the geometry of error amplification.

\section{Proofs}\label{Sec:proofs}

The proof of Theorems \ref{thm:main1} and \ref{thm:main2}, given below, is based on a detailed explicit description of the Prony curves, on one hand, and of a behavior of the moments $m_k, \ k\ge 3$ on these curves, on the other.

\subsection{Parametrization of the Prony curves}\label{sec:parm.Prony.curves}

We denote by $S^x_{F,i}$ the projection of the Prony curve $S_{F,i}$ onto the node plane spanned by the node coordinates $x'_i,x'_{i+1}$.

\bt\label{thm:two.nodes}
The curve $S^x_{F,i}$ is a hyperbola in the plane $x'_i,x'_{i+1}$ defined by the equation
$$
m_0(F)x'_ix'_{i+1}-m_1(F)(x'_i+x'_{i+1})+m_2(F)=0.
$$
The original curve $S_{F,i}$ is parametrized through $x'_i,x'_{i+1}$ in $S^x_{F,i}$ as
$$
a'_i=\frac{m_0(F)x'_{i+1}-m_1(F)}{x'_{i+1}-x'_i}, \ a'_{i+1} = \frac{-m_0(F)x'_i+m_1(F)}{x'_{i+1}-x'_i}, \ (x'_i,x'_{i+1})\in S^x_{F,i}.
$$
\et
\pr
Since all the nodes and amplitudes are fixed on the Prony curve $S_{F,i}$, but $(a_i,x_i),(a_{i+1},x_{i+1})$, we can work only with the partial signals $a_i\delta (x-x_i)+a_{i+1}\delta (x-x_{i+1})$. In other words, we can consider the case of exactly two nodes, i.e. signals of the form $F(x)=a_1\delta (x-x_1)+a_2\delta (x-x_2)\in {\cal P}_2$. In this case there is only one choice $i=1$ for the index $i$ in the definition of the Prony curves $S_{F,i}$, and we denote them by $S_F$.

\smallskip

Alternatively, we can consider algebraic curves $S(m_0,m_1,m_2)$ in ${\cal P}_2$, defined by the equations

\be\label{eq:Prony.fol.2.111}
\begin{array}{c}
a_1+a_2=m_0,\\
a_1x_1+a_2x_2=m_1,\\
a_1x^2_1+a_2x^2_2=m_2,\\
\end{array}
\ee
for any moments $m_0,m_1,m_2$. If we put $m_k=m_k(F), \ k=0,1,2,$ we get $S(m_0,\allowbreak m_1,m_2)=S_F.$ 
Since we are interested in the behavior of the nodes $x_1,x_2$ along the curve $S$, 
we will consider also the node parameter space ${\cal P}^x_2=\{(x_1,x_2)\},$ 
and the projections $S^x(m_0,m_1,m_2) \subset {\cal P}^x_2$ of the Prony curves $S(m_0,m_1,m_2)$ to the node space ${\cal P}^x_2$. 
The following proposition is an extended version of Theorem  \ref{thm:two.nodes}.

\bp\label{Prop:two.nodes11}
The curves $S^x(m_0,m_1,m_2) \subset {\cal P}^x_2$ are hyperbolas in the plane $x_1,x_2$ defined by the equation

\be\label{eq:Prony.fol.2.11}
m_0x_1x_2-m_1(x_1+x_2)+m_2=0.
\ee
They form a two-parametric family, depending only on the ratio of the moments $(m_0:m_1:m_2)$. The corresponding curves $S(m_0,m_1,m_2) \subset {\cal P}_2$ are parametrized as

\be\label{eq:Prony.fol.2.21}
a_1=\frac{m_0x_2-m_1}{x_2-x_1}, \ a_2 = \frac{-m_0x_1+m_1}{x_2-x_1}, \ (x_1,x_2)\in S^x(m_0,m_1,m_2).
\ee
\ep
\pr
We get from the first two equations of (\ref{eq:Prony.fol.2.111}) the following expressions for $a_1,a_2$ through $x_1,x_2$:

\be\label{eq:Prony.fol.121}
a_2=m_0-a_1, \ a_1x_1+(m_0-a_1)x_2=m_1,
\ee
and hence

\be\label{eq:Prony.fol.131}
a_1=\frac{m_0x_2-m_1}{x_2-x_1}, \ a_2 = \frac{-m_0x_1+m_1}{x_2-x_1}.
\ee
The curve $S(m_0,m_1,m_2)$ is defined by all the three equations of (\ref{eq:Prony.fol.2.111}). Substituting (\ref{eq:Prony.fol.131}) into the last equation of (\ref{eq:Prony.fol.2.111}), we see that the projection $S^x= S^x(m_0,m_1,m_2)$ of $S(m_0,m_1,m_2)$ onto the $(x_1,x_2)$-subspace ${\cal P}^x_2 \subset {\cal P}_2$ is obtained in ${\cal P}^x_2$ as the solution of the third degree equation

\be\label{eq:Prony.fol.141}
\frac{m_0x_2-m_1}{x_2-x_1}x^2_1+\frac{-m_0x_1+m_1}{x_2-x_1}x^2_2=m_2.
\ee
An explicit description of the curve $S^x$ can be obtained as follows: we can rewrite the left hand side of equation (\ref{eq:Prony.fol.141}) in the form

\be\label{eq:Prony.fol.151}
\frac{1}{x_2-x_1}[m_0(x^2_1x_2-x_1x^2_2)+m_1(x^2_2-x^2_1)]=-m_0x_1x_2+m_1(x_1+x_2),
\ee
which leads to the equation

\be\label{eq:Prony.fol.161}
m_0x_1x_2-m_1(x_1+x_2)+m_2=0
\ee
for the curve $S^x(m_0,m_1,m_2).$ So this curve is a hyperbola with the center at the point $(\frac{m_1}{m_0},\frac{m_1}{m_0})$, and with the asymptotes $x_1=\frac{m_1}{m_0}, \ x_2=\frac{m_1}{m_0}.$ Equation \ref{eq:Prony.fol.161} is homogeneous in $(m_0,m_1,m_2)$ and hence its solution depends only on the ratio of the moments $(m_0:m_1:m_2)$. Applying expressions \ref{eq:Prony.fol.131} we complete the proof of Proposition \ref{Prop:two.nodes11} and of Theorem \ref{thm:two.nodes}. $\square$ $\square$

\medskip

We expect that the explicit description of the Prony curves given above, can be combined with the analysis of the Prony mapping from the point of view of Singularity Theory, given in \cite{Bat.Yom3,Yom}, including, in particular, repsentation of signals $F$ in the ``bases of finite differences'' introduced in \cite{Bat.Yom3,Yom}.

\subsection{Moments on the Prony curves}

In this section we describe the behavior of the moments $m_k(F), \ k\geq 3,$ along the Prony curve $S_{F,i}$, on its $h$-local part $S_{F,i}(h).$

\bt\label{thm:moments.on.S}
Let $F=\sum_{j=1}^{d}a_{j}\delta\left(x-x_{j}\right)\in {\cal P}_d$ form an $(i,h,M)$-cluster, and let
$\kappa=\frac{x_{i+1}+x_i}{2}$ be the center of the interval $[x_i,x_{i+1}]$. Then for any $F' \in S_{F,i}(h)$ we have
$m_k( F')-m_k(F)=0, \ k=0,1,2,$ while $$ |\ m_k(F')-m_k(F)|\leq 18M(1+|\kappa|)^kh^3, \ k\geq 3.
$$
\et
\pr
As in the previous section, it is sufficient to consider the case of exactly two nodes. By the assumptions, for the signal $F(x)=a_1\delta (x-x_1)+a_2\delta (x-x_2)\in {\cal P}_2$ we have $x_2=x_1+h$, and $|a_1|, |a_2| \leq M$. To simplify the expressions we shall assume that $h\leq 1$. Let us show first that the distance between the nodes remains uniformly bounded from below along $S^x_F(h)$.

\bl\label{lem:bdd.dist11}
For each $F'(x)=a'_1\delta (x-x'_1)+a'_2\delta (x-x'_2)\in S_F(h)$ we have
$$
x'_2-x'_1 \geq \frac{1}{4}h.
$$
\el
\pr
The point $(x_1,x_2)\in {\cal P}^x_2$ is at the distance $\frac{1}{\sqrt 2}h$ from the diagonal $\{x_1=x_2\}$. So the disk $D$ is at the distance $\kappa = (\frac{1}{\sqrt 2}-\frac{1}{2})h > 0.2h$ from the diagonal. Therefore for any $(x'_1,x'_2)\in D$ we have $x'_2-x'_1>0.2\sqrt 2 h > \frac{1}{4}h.$ In particular, this is true for each point of $S^x_F(h)$. $\square$

\smallskip

Next we show that the amplitudes $a'_1, a'_2$ are uniformly bounded on $S_F(h)$.

\bl\label{lem:bdd.ampl11}
For each $F'(x)=a'_1\delta (x-x'_1)+a'_2\delta (x-x'_2)\in S_F(h)$ we have
$$
|a'_1|, |a'_2| \leq 8M.
$$
\el
\pr
By expressions (\ref{eq:Prony.fol.2.21}) in Proposition \ref{Prop:two.nodes11} we have

\be\label{eq:Prony.fol.2.1311}
a'_1=\frac{m_0 x'_2-m_1}{x'_2-x'_1}, \ a'_2 = \frac{-m_0 x'_1+m_1}{x'_2-x'_1}.
\ee
We can write
$$
a'_1=\frac{m_0 x'_2-m_1}{x'_2-x'_1}=\frac{m_0x_2-m_1}{x'_2-x'_1}+\frac{m_0(x'_2-x_2)}{x'_2-x'_1},
$$
and hence
$$
|a'_1|\leq |\frac{m_0x_2-m_1}{x_2-x_1}|\cdot |\frac{x_2-x_1}{x'_2-x'_1}|+|m_0|\cdot |\frac{x'_2-x_2}{x'_2-x'_1}|,
$$
or

\be\label{eq:Prony.fol.2.14}
|a'_1|\leq |a_1|\cdot|\frac{x_2-x_1}{x'_2-x'_1}|+|m_0|\cdot |\frac{x'_2-x_2}{x'_2-x'_1}|.
\ee
By Lemma \ref{lem:bdd.dist11} we have $x'_2-x'_1 \geq \frac{1}{4}h,$ while by the assumptions $x_2-x_1=h$. Since the point $(x'_1,x'_2)$ belongs to the disk $D$ of radius $\frac{1}{2}h$ centered at $(x_1,x_2)$, we have also $|x'_2-x_2|\leq \frac{1}{2}h.$ Therefore (\ref{eq:Prony.fol.2.14}) implies

\be\label{eq:Prony.fol.2.15}
|a'_1|\leq 4|a_1|+2|m_0|\leq 8M,
\ee
since by the assumptions $|a_1|, |a_2| \leq M$, and hence $|m_0|=|a_1+a_2|\leq 2M$. The bound for $|a'_2|$ is obtained exactly in the same way. $\square$

\smallskip

In order to estimate the differences $m_k(F')-m_k(F)$ we now shift the origin into the middle point
$\kappa=\frac{x_1+x_2}{2}$ between the nodes $x_1,x_2$. For $$
F(x)=\sum_{j=1}^d a_j\delta (x-x_j)\in {\cal P}_d
$$
denote by $F^\kappa(x)$ the shifted signal $F^\kappa(x)=F(x-\kappa)$.

\medskip

The following proposition describes the action of the coordinate shift on the moments of general spike-trains (of course, this results remains valid for the moments of any measure on $\mathbb R$).

\bp\label{prop:shift}
$$
m_k(F)=\sum_{l=0}^k \binom{k}{l}(-\kappa)^{k-l} m_l(F^\kappa), \ m_k(F^\kappa) = \sum_{l=0}^k \binom{k}{l}(\kappa)^{k-l} m_l(F).
$$
\ep
\pr
$$
m_k(F^\kappa)= \sum_{j=1}^d a_j (\kappa+x_j)^k =\sum_{j=1}^d a_j \sum_{l=0}^k \binom{k}{l}\kappa^{k-l}x_j^l =
$$
$$
= \sum_{l=0}^k \binom{k}{l}\kappa^{k-l} \sum_{j=1}^d a_j x_j^l=	\sum_{l=0}^k \binom{k}{l}\kappa^{k-l} m_l(F).
$$	
Replacing $\kappa$ by $-\kappa$ we get the second expression. $\square$

\smallskip

Finally we come to estimating the differences $m_k(F')-m_k(F)$. Since by Proposition \ref{prop:shift} the shifted moments are expressed through the original moments of the same and of smaller orders, we see that along the curve $S_F$ the first three shifted moments do not change. Applying Proposition \ref{prop:shift} in the opposite direction, we can write, for $k\geq 3$,

\be\label{eq:shifted.diff}
|m_k(F')-m_k(F)| \leq \sum_{l=3}^k \binom{k}{l}|\kappa|^{k-l}|m_l(F'^\kappa)-m_l(F^\kappa)|.
\ee
By the choice of $\kappa$ we have $|x_1-\kappa|=|x_2-\kappa|=h/2.$ For $(x'_1,x'_2)\in D$ we have $|x'_1-\kappa|,|x'_2-\kappa|\leq h$. Hence we
obtain, using Lemma \ref{lem:bdd.ampl11}, $$
|m_l(F^\kappa)|=|a_1(x_1-\kappa)^l+a_2(x_2-\kappa)^l|\leq 2M(\frac{h}{2})^l,
$$
$$
|m_l(F'^\kappa)|=|a'_1 (x'_1-\kappa)^l+a'_2(x'_2-\kappa)^l|\leq 16Mh^l.
$$
Consequently, $|m_l(F'^\kappa)-m_l(F^\kappa)|\leq (16+2(\frac{1}{2})^l)Mh^l \leq 18Mh^l.$ Substituting this into equation
(\ref{eq:shifted.diff}) we get

$$
|m_k(F')-m_k(F)| \leq 18M\sum_{l=3}^k \binom{k}{l}|\kappa|^{k-l}h^l \leq
$$
$$
\leq 18Mh^3\sum_{l=3}^k \binom{k}{l}|\kappa|^{k-l}\leq 18M(1+|\kappa|)^kh^3.
$$
This completes the proof of Theorem \ref{thm:moments.on.S}. $\square$

\subsection{Proof of Theorem \ref{thm:main1}}

We have to show that for each $\e\ge C(F)h^3$, with $C(F)=18M(1+|\kappa|)^N,$ the $\e$-error set $E_\e(F)$ contains the local
Prony curve $S_{F,i}(h).$ By Theorem \ref{thm:moments.on.S} we have for any $F'\in S_{F,i}(h)$ and for each $k\leq N$ $$
|\ m_k(F')-m_k(F)|\leq 18M(1+|\kappa|)^kh^3 \leq 18M(1+|\kappa|)^Nh^3=C(F)h^3 \le \e,
$$
and therefore $S_{F,i}(h)\subset E_\e(F).$ This completes the proof. $\square$

\subsection{Proof of Theorem \ref{thm:main2}}

By definition, for $F=(a,x)\in {\cal P}_d$ the worst case error $\rho^x_\e(F)$ in reconstruction of the nodes of $F$ is
$$
\rho^x_\e(F)=\sup_{F'=(a',x')\in E_\e(F)} ||x'-x||.
$$
The projection $S^x_{F,i}(h)$ of the $h$-local Prony curve $S_{F,i}(h)$ to the coordinate plane of $(x'_i,x'_{i+1})$ is a hyperbola, passing through the point $(x_i,x_{i+1})$, and it crosses the boundary of the disk $D$ of radius $\frac{1}{2}h$ centered at $(x_i,x_{i+1})$, at exactly two points. Let $F''=(a'',x'')$ be one of the corresponding endpoints of $S_{F,i}(h).$ Then the distance between the nodes of $F''$ and the nodes of $F$ is exactly $\frac{1}{2}h$. By Theorem \ref{thm:main1} we have $S_{F,i}(h)\subset E_\e(F),$ and therefore $F''\in E_\e(F)$. We conclude that
$$
\rho^x_\e(F)=\sup_{F'=(a',x')\in E_\e(F)} ||x'-x||\ge ||x''-x||=\frac{1}{2}h.
$$
This completes the proof. $\square$






\section{A case of complex nodes}\label{sec:compl.nodes}

The goal of this section is to compare two approaches to the reconstruction problem for real spike-train signals with exactly two nodes: from their moments, and from their Fourier samples, and to pose some related open questions. For a signal
$$
F(x)=a_1\delta(x-x_1)+a_2\delta(x-x_2)\in {\cal P}_2
$$
we have for its Fourier transform $f_s(F):={\cal F}(F)(s)=a_1e^{isx_1}+a_2e^{isx_2}$. Taking samples $f_k(F)$ at the points $s=0,1,\ldots,k,\ldots,$ we get $f_k(F)=a_1e^{ikx_1}+a_2e^{ikx_2}$.

We see immediately, that the Fourier samples  $f_k(F)=a_1e^{ikx_1}+a_2e^{ikx_2}$ coincide with the moments $m_k(\tilde F)$ for a signal $\tilde F(x)=a_1\delta(x-e^{ix_1})+a_2\delta(x-e^{ix_2})$ with the {\it complex nodes} $e^{ix_1},e^{ix_2}$.

Recently in \cite{Aki.Gol.Yom} a trigonometric reconstruction method for $\tilde F$ from their moments $m_k(\tilde F)$ was introduced, which uses the following four real measurements:
$$
|m_0|,|m_1|,|m_2|, \ \text {and the imaginary part} \ \Im m_1.
$$
Notice that each complex moment provides (at least, formally) {\it two real measurements}: its real and its complex parts. So taking four moments $m_k, \ k=0,1,2,3,$ as in a true complex Prony system, gives us eight real equations, while the signals $F$ and $\tilde F$ have only four real parameters: $a_1,a_2,x_1,x_2$.

\smallskip

This leads us to the following question: {\it what real measurements (coming from the real or from the complex moments) do we really need? Can we improve the reconstruction accuracy by a ``correct choice'' of the measurements?}

\smallskip

The last question is directly connected to the main results of the current paper, because of the following fact: the trigonometric reconstruction method of \cite{Aki.Gol.Yom} uses as an input only {\it three complex moments $m_0,m_1,m_2$}. According to the approach of the present paper, we would expect for $F$ (or $\tilde F)$ with two nodes in a distance $h\ll 1$ the worst case node error amplification factor to be of order $(\frac{1}{h})^{l-2}$, where $l$ is the number of the moments used. For the trigonometric method $l=3$, and this would lead to the amplification factor of order $\frac{1}{h}$, while for the Prony inversion it is shown above to be of order $(\frac{1}{h})^{2}$ - an apparent contradiction.

Our initial experiments (partially reported in \cite{Aki.Gol.Yom}) indicate also for the trigonometric method that the worst case error amplification factor of order $(\frac{1}{h})^{2}$. This leads to the following question, which may be important for better understanding the patterns of error amplification in different methods of spike-train reconstruction:

\smallskip

\noindent {\it Is it possible to extend the approach of the present paper to the analysis of the error amplification in trigonometric reconstruction? Where do we (presumably) lose, in the trigonometric method, the accuracy gained by not using the fourth moment?}

\smallskip

For the reader convenience we shortly recall below the main steps of the trigonometric reconstruction method of \cite{Aki.Gol.Yom}.

\subsection{Trigonometric reconstruction: main steps}\label{sec:trig.reconstr}

Consider signal $F(u)$ with {\it complex} variables $(x, y)$ of the form (in this section, we use notations from \cite{Aki.Gol.Yom}):

\be\label{eq:unit.circle}
F(u)=a\delta(u-x)+b\delta(u-y), \ x=e^{i\phi}, \ y=e^{i\theta}, \ \phi,\theta,a,b \in {\mathbb R}
\ee
Our measurements are the complex moments

\be\label{eq:unit.circle2}
m_k=\int_{\mathbb C} x^k F(x)dx=ae^{ik\phi}+be^{ik\theta}.
\ee
We can write the moments for signal \eqref{eq:unit.circle} in the following form:
\be\label{eq:unit.circle3}
  m_k = ax^k + by^k.
\ee

\subsection{Recovery of phase difference}

Introduce the following definitions for phases of $x$ and $y$:
$$
  \phi = -2\pi \mu, \quad
  \theta = -2\pi \nu, \quad
  \Delta = \phi - \theta.
$$

From \eqref{eq:unit.circle}, we have
$$
  x = e^{i \phi} = \cos \phi + i \sin \phi, \quad
  y = e^{i \theta}= \cos \theta + i \sin \theta.
$$

Consider real and imaginary parts of the moment $m_k$:

\begin{equation}
  \label{eq:TProny2D-ReIm}
  \begin{cases}
    \Re m_k = a \cos k\phi + b \cos k\theta \\
    \Im m_k = a \sin k\phi + b \sin k\theta
  \end{cases}
\end{equation}

Now, we get

$$
  \begin{cases}
  (\Re m_k)^2 = a^2 \cos^2 k\phi + b^2 \cos^2 k\theta + 2ab \cos k\phi \cos k\theta \\
  (\Im m_k)^2 = a^2 \sin^2 k\phi + b^2 \sin^2 k\theta + 2ab \sin k\phi \sin k\theta
  \end{cases}
$$

Let $M_k = |m_k|$. So

$$
  M_k^2 = |m_k|^2 = a^2 + b^2 + 2ab \cos k\Delta.
$$

Since $M_0^2 = a^2 + 2ab + b^2$, we get

\begin{equation}
  \label{eq:temp1}
  2\sin^2 \frac{k\Delta}{2}= -\dfrac{M_k^2 - M_0^2}{2ab}.
\end{equation}

It follows from \eqref{eq:temp1} that

$$
  \dfrac{\sin^2 \Delta/2}{\sin^2 \Delta} \equiv
  \dfrac{1}{4 \cos^2(\Delta/2)}
  = \dfrac{M_1^2 - M_0^2}{M_2^2 - M_0^2}.
$$

and hence,
$$
  \dfrac{1+\cos \Delta}{2} =
  \dfrac{1}{4} \dfrac{M_2^2 - M_0^2}{M_1^2 - M_0^2}
  \quad \text{or} \quad
  \Delta = \arccos \Bigg(
    \dfrac{2 M_1^2 - M_0^2 - M_2^2}{2(M_0^2 - M_1^2)}.
  \Bigg)
$$

\subsection{Amplitudes recovery}

Recall that

$$
  \cos \Delta - 1= \dfrac{M_1^2 - M_0^2}{2ab} \qquad
  \text{and} \qquad
  \cos\Delta + 1 =
  \dfrac{1}{2} \dfrac{M_2^2 - M_0^2}{M_1^2 - M_0^2}.
$$

Thus

$$
  2 = \dfrac{1}{2}
  \dfrac{M_2^2 - M_0^2}{M_1^2 - M_0^2} - \dfrac{M_1^2 - M_0^2}{2ab}
  \quad \text{or} \quad
  \dfrac{M_1^2 - M_0^2}{ab} =
  \dfrac{M_2^2 + 3M_0^2 - 4M_1^2}{M_1^2 - M_0^2}.
$$

So we get

$$
  ab = \dfrac{(M_1^2 - M_0^2)^2}{M_2^2 + 3M_0^2 - 4M_1^2}.
$$

Now, in order to find the unknown amplitudes, we have to solve the system

$$
  \begin{cases}
    a + b = M_0, \\
    ab = \dfrac{(M_1^2 - M_0^2)^2}{M_2^2 + 3M_0^2 - 4M_1^2},
  \end{cases}
$$

which is reduced to the quadratic equation

$$
  a^2 - M_0 a +
  \dfrac{(M_1^2 - M_0^2)^2}{M_2^2 + 3M_0^2 - 4M_1^2} = 0
$$

with the discriminant

$$
  D = M_0^2 - 4\dfrac{(M_1^2 - M_0^2)^2}{M_2^2 + 3M_0^2 - 4M_1^2}
$$

Now, we obtain the amplitudes:

$$
  a = \dfrac{M_0 \pm \sqrt{D}}{2}, \quad
  b = \dfrac{M_0 \mp \sqrt{D}}{2}.
$$

\subsection{Phases recovery}

Let $\phi = \theta + \Delta$. Then the system \eqref{eq:TProny2D-ReIm} has the form

$$
  \begin{cases}
    \Re m_1 = a \cos (\theta + \Delta) + b \cos \theta, \\
    \Im m_1 = a \sin (\theta + \Delta) + b \sin \theta
  \end{cases}
$$

or

$$
  \begin{cases}
    \Re m_1 = a \Big( \cos\theta \cos\Delta - \sin\theta \sin\Delta \Big) + b \cos \theta, \\
    \Im m_1 = a \Big( \sin\theta \cos\Delta + \cos\theta \sin\Delta \Big) + b \sin \theta
  \end{cases}
$$

and

$$
  \begin{cases}
    \Re m_1 = (a\cos\Delta + b)\cos\theta + (-a\sin\Delta) \sin\theta, \\
    \Im m_1 = (a\sin\Delta) \cos\theta + (a\cos\Delta + b)\sin\theta.
  \end{cases}
$$

Thus, we get

$$
  \theta = -\arccos \Big(
    \dfrac
      {\Re m_1 (a\cos\Delta + b) + \Im m_1 (a\sin\Delta)}
      {(a\cos\Delta + b)^2 + (a\sin\Delta)^2}
  \Big), \quad
  \phi = \theta + \Delta.
$$









\clearpage

\bibliographystyle{amsplain}

\begin{thebibliography}{10}

\bibitem{Aki.Bat.Yom} Akinshin, A., Batenkov, D. Yomdin, Y. {\it Accuracy of spike-train Fourier reconstruction for colliding nodes} in Sampling Theory and Applications (SampTA), Pp. 617--621 (IEEE, 2015). doi:10.1109/SAMPTA.2015.7148965

\bibitem{Aki.Gol.Yom}  A.A.Akinshin, V.P.Golubyatnikov, Y.N.Yomdin {\it Low-dimensional Prony systems} (In Russian), Proc. International Conference ``Lomonosov readings in Altai: fundamental problems of science and education'', Barnaul, 20 -- 24 October 2015, Altai state university. p. 443 - 450.

\bibitem{Aza.Cas.Gam} Jean-Marc Aza{\"i}s, Yohann de~Castro, and Fabrice Gamboa, {\it Spike detection from inaccurate samplings},
Applied and Computational Harmonic Analysis, in press.

\bibitem{Bat1} D. Batenkov, {\it Accurate solution of near-colliding {Prony} systems via decimation and homotopy continuation},
{arXiv}:1501.00160 [cs, math], December 2014.

\bibitem{Bat.Yom2} D. Batenkov and Y. Yomdin, {\it On the accuracy of solving confluent {Prony} systems},
SIAM J.Appl.Math., 73(1):134--154, 2013.

\bibitem{Bat.Yom3} D. Batenkov and Y. Yomdin, {\it {Geometry} and {Singularities} of the {Prony} mapping},
Journal of Singularities, 10:1--25, 2014.

\bibitem{Bey.Mon2} G. Beylkin and L. Monzon, {\it Nonolinear inversion of a band-limited Fourier transform}, Appl. Comput.
Harmon. Anal. 27, 351--366, 2009.

\bibitem{Blu.Dra.Vet.Mar.Cou} T. Blu, P.-L. Dragotti, M. Vetterli, P. Marziliano, L. Coulot, {\it Sparse Sampling of Signal
Innovations}, IEEE Signal Proc. Magazine, 31--40, March 2008.


\bibitem{Can.Fer1} Emmanuel J. Cand{\`e}s and Carlos Fernandez-Granda, {\it {Super}-{Resolution} from {Noisy} {Data}},
Journal of Fourier Analysis and Applications, 19(6):1229--1254, December 2013.

\bibitem{Can.Fer2} Emmanuel J. Cand{\`e}s and Carlos Fernandez-Granda, {\it {Towards} a {Mathematical} {Theory} of {Super}-resolution},
Communications on Pure and Applied Mathematics, 67(6):906--956, June 2014.

\bibitem{Con.Hir} L. Condat, A. Hirabayashi, {\it A New Projection Method for the Recovery of Dirac Pulses from Noisy Linear
Measurements}, https://hal.archives-ouvertes.fr/hal-00759253v5, 2014.

\bibitem{Dem.Ngu} L. Demanet and N. Nguyen, {\it The recoverability limit for superresolution via sparsity}, Preprint, 2014.

\bibitem{Dem.Nee.Ngu} Laurent Demanet, Deanna Needell, and Nam Nguyen, {\it Super-resolution via superset selection and pruning},
In Proceedings of the 10th International Conference on Sampling Theory and Applications ({SAMPTA}), 2013.

\bibitem{Don} D.L. Donoho, {\it Superresolution via sparsity constraints}, {SIAM} Journal on Mathematical Analysis, 23(5):1309--1331,
1992.

\bibitem{Duv.Pey} Vincent Duval and Gabriel Peyr{\'e}, {\it Exact support recovery for sparse spikes deconvolution},
{arXiv} preprint {arXiv}:1306.6909, 2013.

\bibitem{Fer} Carlos Fernandez-Granda, {\it Support detection in super-resolution}, In Proc. of 10th Sampling Theory and Applications
({SAMPTA}), pages 145--148, 2013.

\bibitem{Hec.Mor.Sol} Reinhard Heckel, Veniamin I. Morgenshtern, and Mahdi Soltanolkotabi, {\it {Super}-{Resolution} {Radar}},
{arXiv}:1411.6272 [cs, math], November 2014.

\bibitem{Lev.Ful} S. Levy and P.K. Fullagar, {\it Reconstruction of a sparse spike train from a portion of its spectrum
and application to high-resolution deconvolution}, Geophysics, 46(9):1235--1243, 1981.

\bibitem{Lia.Fan} Wenjing Liao and Albert Fannjiang, {\it {MUSIC} for {Single}-{Snapshot} {Spectral} {Estimation}: {Stability}
and {Super}-resolution}, {arXiv}:1404.1484 [cs, math], April 2014.

\bibitem{Mor.Can} V. I. Morgenshtern and E. J. Candes, {\it Stable super-resolution of positive sources: the discrete setup},
Preprint, 2014.


\bibitem{Moi} Ankur Moitra, {\it The {Threshold} for {Super}-resolution via {Extremal} {Functions}}, {arXiv}:1408.1681
[cs, math, stat], August 2014.

\bibitem{Ode.Bar.Pis} J. Odendaal, E. Barnard, and C. W. I. Pistorius, {\it Two-dimensional superresolution radar imaging using the
{MUSIC} algorithm}, IEEE Transactions on Antennas and Propagation, 42(10):1386--1391, 1994.

\bibitem{Pet.Plo} T. Peter, G. Plonka, {\it A generalized Prony method for reconstruction of sparse sums of
eigenfunctions of linear operators}, Inverse Problems 29 (2013), 025001.

\bibitem{Pet.Pot.Tas} T.~Peter, D.~Potts, and M.~Tasche.
\newblock Nonlinear approximation by sums of exponentials and translates.
\newblock {\em SIAM Journal on Scientific Computing}, 33(4):1920, 2011.

\bibitem{Plo.Wis} Gerlind Plonka, M. Wischerhoff, {\it How many Fourier samples are needed for real function reconstruction?},
J. Appl. Math. Comput. 42, 117--137, 2013.

\bibitem{Pot.Tas} Daniel Potts and Manfred Tasche, {\it Fast ESPRIT algorithms based on partial singular value decompositions}


\bibitem{Yom} Y. Yomdin, {\it Singularities in algebraic data acquisition}, Real and complex singularities,
London Math. Soc. Lecture Note Ser., 380:378--396, 2010.

\bibitem{Aki} A. A. Akinshin {\it Prony analysis [Technical report]}, authorea.com/51974


\end{thebibliography}

\end{document}